\documentclass[11pt,a4paper,leqno]{article}
\usepackage{a4wide}
\usepackage{latexsym}
\usepackage{amsmath,amsthm}
\usepackage{amssymb}
\usepackage{stackrel}  
\usepackage{color}
\usepackage{graphicx}
\usepackage[linesnumbered,ruled,vlined]{algorithm2e}

\usepackage{etoolbox}
\usepackage{hyperref}

\makeatletter
\newcounter{author}
\renewcommand*\author[1]{%
  \stepcounter{author}%
  \ifnum\c@author=1
    \gdef\@author{#1}%
  \else
    \xdef\@author{\unexpanded\expandafter{\@author\and#1}}%
  \fi
  \csgdef{author@\the\c@author}{#1}}
\newcommand*\email[1]{%
  \csgdef{email@\the\c@author}{#1}}
\newcommand*\address[1]{%
  \csgdef{address@\the\c@author}{#1}}
\AtEndDocument{%
  \xdef\author@count{\the\c@author}%
  \c@author=1
  \print@authors}
\newcommand*\print@authors{%
  \ifnum\c@author>\author@count
  \else
    \print@author{\the\c@author}%
    \advance\c@author by 1
    \expandafter\print@authors
  \fi}
\newcommand*\print@author[1]{%
  \par\medskip
  \begin{tabular}{@{}l@{}}%
    \textsc{\textbf{\csuse{author@#1}}}\\
    \csuse{address@#1}\\
    \textit{E-mail}:
    \href{mailto:\csuse{email@#1}}{\csuse{email@#1}}
  \end{tabular}}
\makeatother

\newtheorem{lemma}{Lemma}

\newtheorem{theo}{Theorem}
\newtheorem{corol}{Corollary}
\theoremstyle{definition}
\newtheorem{example}{Example}
\newtheorem{remark}{Remark}
\newcommand{\C}{\mathbb{C}}
\newcommand{\N}{\mathbb{N}}

\newcommand{\R}{\mathbb{R}}

\newcommand{\Oo}{\mathcal{O}}

\newcommand{\ord}{\mathrm{ord}}
\newcommand{\conv}{\mathrm{conv}}

\newcommand{\corner}{\mathbin{\vrule height 0.13ex depth 0pt width 1.3ex
\vrule height 1.6ex depth 0pt width
0.13ex}\,}

\title{Formal power series solutions with coefficients defined on shrinking discs for some partial differential equations}

\author{Alberto Lastra}
\address{Universidad de Alcal\'a, Departamento de F\'isica y Matem\'aticas, E-28871.\\ 
Alcal\'a de Henares, Madrid, Spain}
\email{alberto.lastra@uah.es}

\author{S{\l}awomir Michalik}
\address{Faculty of Mathematics and Natural Sciences, College of Science, Cardinal Stefan\\ Wyszy{\'n}ski University in Warsaw, W\'oycickiego 1/3, 01-938 Warszawa, Poland;\\
Shibaura Institute of Technology,
Department of Engineering and Design,\\
Saitama 337-8570, Japan}
\email{s.michalik@uksw.edu.pl}

\author{Maria Suwi{\'n}ska}
\address{Faculty of Mathematics and Natural Sciences, College of Science, Cardinal Stefan \\Wyszy{\'n}ski  University in Warsaw, W\'oycickiego 1/3, 01-938 Warszawa, Poland}
\email{m.suwinska@op.pl}

\date{}

\begin{document}
\maketitle

\thispagestyle{empty}
{ \small \begin{center}
{\bf Abstract}
In this paper conditions, under which an integro-differential operator is a linear automorphism, are provided. Alternatively, the problem can be considered in terms of existence of a unique formal power series solution for a linear Cauchy problem, where the coefficients of the solution are of a certain Gevrey order on progressively shrinking domains.
\end{center}

\smallskip

}
\bigskip \bigskip

\section{Introduction}

In the present work we study the integro-differential operators of the form $P(\partial_t,\partial_z)\partial_t^{-m}$ with $P$~standing for a polynomial with coefficients from $\C[[t,z]]$. The problem addressed in this paper concerns the subspace $(G_s(r(n))_{n\in\N_0}[[t]]$ of $\C[[z]]_s[[t]]$ with coefficients that are of a fixed Gevrey order $s\ge 0$, but in progressively shrinking domains. In particular, we aim to give conditions, under which $P(\partial_t,\partial_z)\partial_t^{-m}$ can be extended to a linear automorphism on $(G_s(a(n+1)^{-\alpha}))_{n\in\N_0}[[t]]$ for fixed $a>0$ and $\alpha>0$. 

The study of formal power series, whose coefficients are analytic (or, more generally, of a certain Gevrey order) in shrinking domains, is very recent. However, it quickly gained much interest from researchers (see \cite{cala,lama4,tahara3}). Such problems naturally occur while studying $q$-difference and $q$-difference differential equations and their formal solutions, i.e., solutions given in the form of a power series that can potentially be convergent everywhere except for the origin. In the present paper we aim to extend this natural notion from the $q$-difference differential equations towards classical partial differential operators.

 Moreover, it is also worth noting that determining existence and uniqueness of formal power series solutions for partial differential equations is significant in itself. Particularly, it is an important first step in the procedure of summability, where an analytic solution can be derived by means of applying Borel and Laplace transforms to the formal solution (see for example \cite{ba2,lamisu,mi,Remy1,remy22} and references therein). In this process, knowledge of the Gevrey upper bounds for coefficients of the formal solution is the key to determining the appropriate Borel transform. Such results are known as Maillet type theorems and are also a relevant topic of interest. For more information we refer the reader to, among others, \cite{LastraTahara,shirai15}.

The paper is structured in the following way. First, we establish the notation and recall some relevant definitions in Section \ref{section:notation}.  
Section \ref{section:3} is focused on a linear differential operators of the form
\begin{equation*}
	P(n,\partial_z):=\sum_{(i,j)\in\Lambda}w_{ij}(n)a_{ij}(z)z^j(z\partial_z)^i,
\end{equation*}
where $w_{ij}$ are polynomials, $a_{ij}(z)$ are of Gevrey order $s$ on a disc at the center in the origin with fixed radius $r>0$, and $a_{ij}(0)\neq 0$. In Theorem \ref{th:1} we give estimates for coefficients of a formal power series solution of equation
\begin{equation*}
	P(n,\partial_z)u_n(z)=f_n(z)\quad\text{for every}\quad n\in\N_0,
\end{equation*}
under certain conditions for the operator $P$ and the inhomogeneity $f_n$.

The paper will be concluded by Section \ref{section:4}, where Theorem \ref{th:1} will be combined with previous results from \cite{lamisu5} and applied towards integro-differential operators of the form
\begin{equation*}
	P(\partial_t,\partial_z)=\sum_{(q,r)\in\Delta}a_{qr}(t,z)\partial_t^q\partial_z^r,
\end{equation*}
where $\Delta\subset\N_0\times\N_0$ is a finite set of indices, all $a_{qr}(t,z)$ are given by formal power series in $t$ with coefficients from $G_s(r)$. Assuming that $m:=\max_{(q,r)\in\Gamma}(q-\ord_t a_{qr})\ge 0$, where $\ord_t a_{qr}$ denotes the order of zero at $t=0$ for each $a_{qr}$, in Theorem \ref{th:4} we give the conditions, under which the operator $P(\partial_t,\partial_z)\partial_t^{-m}$ is a linear automorphism on $\C[[t,z]]$ and can be extended on $(G_s(rn^{-\alpha}))_{n\in\N}[[t]]$ for certain $\alpha>0$. We achieve this by reducing the problem to analyzing operators of the form $P(n,\partial_z)$ from the previous section. We finish the paper by adding some final comments, remarks and potential directions for further research in Section \ref{section:final_remarks}. For the sake of clarity, all auxiliary lemmas used in the main parts of the paper are gathered in Section \ref{section:lemmas}.

\section{Notation and preliminary definitions}\label{section:notation}
Let us now introduce the notation used throughout the paper and recall some necessary definitions related to the concept of a Gevrey order.

By $\N$ we denote the set of all positive integers. Then, $\N_0=\N\cup\{0\}$.

Let $r>0$. We call the set $D_r=\{z\in\C\ |z|<r\}$ an open disc with a center at the origin and the radius $r$.

The set of all formal power series in variable $t$ with coefficients from a given space $\mathbb{E}$ is denoted by $\mathbb{E}[[t]]$.

Let us consider a formal power series $\sum_{n=0}^\infty u_n z^n\in\C[[z]]$. If there exist positive constants $A,C>0$ such that 
$$
|u_n|\le AC^n\Gamma(1+sn)\textrm{ for every }n\in\N_0,
$$
where $\Gamma(\cdot)$ denotes the gamma function and $s\geq 0$, then we say that this formal power series is of Gevrey order $s$. The space of all formal power series of order $s$ is denoted by $\C[[z]]_s$.

For given $s\geq 0$ and $r>0$ we will use the following subspace of $\C[[z]]_s$:
\begin{equation*}
 G_s(r):=\Big\{\sum_{n=0}^{\infty}u_n z^n\in\C[[z]]_s \colon \sum_{n=0}^{\infty}\frac{u_n}{\Gamma(1+sn)}z^n\in\Oo(D_r) \Big\},
\end{equation*}
where $\Oo(D_r)$ denotes the space of holomorphic functions on the disc $D_r$. Observe that in the special case $s=0$ the space $G_0(r)$ coincides with $\Oo(D_r)$.

Next, for a given sequence of positive radii $(r(n))_{n\geq 0}$ and for $s\geq 0$, we define the following subspace of $\C[[z]]_s[[t]]$:
\begin{equation*}
 (G_s(r(n)))_{n\in\N_0}[[t]]:=\Big\{\sum_{n=0}^{\infty}u_n(z)t^n\in\C[[z]]_s[[t]] \colon u_n(z)\in G_s(r(n))\ \text{for every}\ n\in\N_0 \Big\}.
\end{equation*}

Lastly, let us recall the definition of the Newton polygon (see \cite{ramis84}). For any point $(a,b)\in\R^2$ we define a set
$$
\corner(a,b)=\{(x,y)\in\R^2:\ x\le a,\,y\ge b\}.
$$
Let $\Lambda\subseteq\N_0\times\N_0$ be a finite set of indices. We also assume that $w_{ij}(n)$ is a polynomial, $a_{ij}(z)\in\C[[z]]$ and $a_{ij}(0)\neq 0$ for every $(i,j)\in\Lambda$. Then the Newton polygon for an operator 
$$
		P(n,\partial_z)=\sum_{(i,j)\in\Lambda}w_{ij}(n)a_{ij}(z)z^j(z\partial_z)^i
$$
is defined as a set
$$
N(P(n,\partial_z))=\conv\left\{\bigcup_{(i,j)\in\Lambda}\corner(i,j)\right\}.
$$
Here, $\conv\{A\}$ denotes the convex hull of any given set $A$.

\section{Differential operators in one variable}\label{section:3}
In this section we will consider a sequence of differential operators in one variable $P(n,\partial_z)$ for $n\in\N_0$. To introduce it
we assume that
$\Lambda\subseteq \N_0\times\N_0$ is a finite set of indices such that its subset $\{i\in\N_0\colon (i,0)\in\Lambda\}$ is nonempty, $p:=\max\{i\in\N_0\colon (i,0)\in\Lambda\}$, $\Lambda':=\{(i,j)\in\Lambda\colon j>0\}$,
\begin{equation*}
s:=\max\left(0,\max_{(i,j)\in\Lambda'}\frac{i-p}{j}\right),
\end{equation*}
$\Lambda_j:=\{i\in\N_0\colon (i,j)\in\Lambda\}$ for every $j\in\N_0$, and $\tilde{\Lambda}_s:=\{(i,j)\in\Lambda\colon js=i-p\}$.

For every $n\in\N_0$ we define a linear differential operator
\begin{equation}\label{eq:P_n}
 P(n,\partial_z):=\sum_{(i,j)\in\Lambda}w_{ij}(n)a_{ij}(z)z^j(z\partial_z)^i,
\end{equation}
where $w_{ij}(n)$ are polynomials, $a_{ij}(z)\in G_s(D_r)$ for some fixed $r>0$, and $a_{ij}(0)\neq 0$ for all $(i,j)\in \Lambda$.

We assume that the numbers
\begin{equation*}
 W(n,k):=\sum_{i\in\Lambda_0}w_{i0}(n)a_{i0}(0)k^i
\end{equation*}
satisfy the following conditions:
\begin{itemize}
\item non-resonance condition:
\begin{equation}\label{eq:nrc}
 W(n,k)\neq0\ \text{for every}\ n,k\in\N_0;
\end{equation}
\item strong non-resonance condition:
\begin{equation}\label{eq:snrc}
 \text{there exists}\ C_0>0\ \text{such that}\ |W(n,k)|\geq C_0\ \text{for every}\ n,k\in\N_0.
\end{equation}
\end{itemize}

\begin{remark}
 Of course the strong non-resonance condition (\ref{eq:snrc}) implies the non-resonance condition (\ref{eq:nrc}), so it is sufficient to assume that the numbers $W(n,k)$ satisfy (\ref{eq:snrc}). 
 
 The converse does not hold for every polynomial $W(x,y)$. As an example it is enough to take a polynomial $W(x,y)=x-\lambda (y+1)$, where $\lambda$ is a fixed positive Liouville number. Since $\lambda$ is an irrational number, the non-resonance condition (\ref{eq:nrc}) holds.
 On the other hand, by the definition of the Liouville number, for every $m\in\N$ there exist $n,k\in\N$ such that
 \begin{equation}\label{eq:liouville}
  \Big|\lambda-\frac{n}{k+1}\Big| < \frac{1}{(k+1)^m}\quad\text{or equivalently}\quad \Big|W(n,k)\Big| < \frac{1}{(k+1)^{m-1}}.
 \end{equation}
 Hence for every $C_0>0$ we may find sufficiently large $m\in\N$ and $n,k\in\N$ satisfying 
 $$|W(n,k)|< \frac{1}{(k+1)^{m-1}}<C_0,$$
 which means that the strong non-resonance condition (\ref{eq:snrc}) does not hold.
\end{remark}

Observe that every operator $P(n,\partial_z)$ given by (\ref{eq:P_n}) may be written in such a form that 
\begin{equation}\label{eq:cond_a_i0}
a_{i0}(z)\equiv a_{i0}(0)\equiv\text{const}\quad \text{for every}\quad i\in\Lambda_0.
\end{equation}
Hence, without loss of generality we may assume that the condition (\ref{eq:cond_a_i0}) holds.


Let us introduce the nonnegative rational number $\alpha$ dependent on the sequence of operators $P(n,\partial_z)$, $n\in\N_0$, given by (\ref{eq:P_n}) in the following way
\begin{equation}\label{eq:alpha}
 \alpha:=\max\left(0, \max_{(i,j)\in\tilde{\Lambda}_s\cap\Lambda'}\frac{\deg w_{ij}-\deg w_{p0}}{j}\right),
\end{equation}
where $\deg w$ denotes the degree of a polynomial $w(n)$.

Now, we are ready to formulate the first main result of the paper.
\begin{theo}\label{th:1}
 We assume that $u_n(z)=\sum_{k=0}^{\infty}u_{n,k}z^k$ is a formal power series solution of the equation
 \begin{equation}\label{eq:n}
  P(n,\partial_z)u_n(z)=f_n(z)\quad\text{for every}\quad n\in\N,
 \end{equation}
 where the inhomogeneity $f_n(z)=\sum_{k=0}^{\infty}f_{n,k}z^k$ is a given formal power series of Gevrey order $s$ and the family of operators $P(n, \partial_z)$ defined by (\ref{eq:P_n}) satisfies (\ref{eq:snrc}).
 
 Then there exists a positive constant $B>0$ depending only on the sequence of operators $(P(n,\partial_z))_{n\in\N_0}$, such that for every sequence of positive numbers $(\tilde{A}(n))_{n\in\N}$  the following implication holds:\\
 if
 \begin{equation}\label{eq:f_nk_est}
  |f_{n,k}|\leq \tilde{A}(n)B^k n^{\alpha k}k!^s\quad\text{for every}\quad k,n\in\N
 \end{equation}
 then
 \begin{equation}\label{eq:u_nk_est}
  |u_{n,k}|\leq A(n)B^k n^{\alpha k}k!^s\quad\text{for every}\quad k,n\in\N
\end{equation}
for some sufficiently large sequence of positive numbers $(A(n))_{n\in\N}$ depending on $(\tilde{A}(n))_{n\in\N}$.

Moreover, if the exponent $\alpha$ defined by (\ref{eq:alpha}) is positive then it is the smallest possible one for which the above implication holds.
\end{theo}

\begin{proof}
We will divide the proof into six steps.

\emph{Step I. Estimation of polynomials $w_{ij}(n)$ and Taylor coefficients of functions $a_{ij}(z)$ from above.}

Since the set $\Lambda'$ is finite, there exists a constant $C_1>0$ such that 
\begin{equation}\label{eq:est_w_ij}
 |w_{ij}(n)|\leq C_1 n^{\deg w_{ij}}\quad \text{for every}\quad (i,j)\in\Lambda'\quad\text{and for every}\quad n\in\N. 
\end{equation}
For the same reason we can estimate the coefficients $a_{ijl}$ of the functions
\begin{equation*}
 a_{ij}(z)=\sum_{l=0}^{\infty}a_{ijl}z^l\in G_s(D_r),\quad\text{where}\quad (i,j)\in\Lambda'.
\end{equation*}
Namely, there exist constants $A_1,B_1>0$ such that
\begin{equation}\label{eq:est_a_ijl}
 |a_{ijl}|\leq A_1 B_1^l (l!)^s\quad\text{for every}\quad (i,j)\in\Lambda'\quad\text{and for every}\quad l\in\N_0.
\end{equation}
\medskip\par
\emph{Step II. The estimation of the polynomials $W(n,k)$ from below.}

First, observe that the growth of $|W(n,k)|$ for $n,k\in\N$ is determined by the terms with maximal powers with respect to $n$, i.e., by the terms of the form
$n^{\deg w_{i0}}k^i$ for $i\in\Lambda_0$. In particular, the leading term is given by $n^{\deg w_{p0}}k^p$ if
\begin{equation*}
 n^{\deg w_{i0}}k^i\leq n^{\deg w_{p0}}k^p \quad\text{for every}\quad i\in\Lambda_0.
\end{equation*}
It means that
\begin{equation*}
 k^{p-i}\geq n^{\deg w_{i0}-\deg w_{p0}}\quad\text{for every}\quad i\in\Lambda_0,
\end{equation*}
or equivalently
\begin{equation*}
 k\geq n^{\frac{\deg w_{i0}-\deg w_{p0}}{p-i}}\quad\text{for every}\quad i\in\Lambda_0\setminus\{p\}.
\end{equation*}
Hence, for
\begin{equation*}
 \gamma:=\max\left(0,\max_{i\in\Lambda_0\setminus\{p\}}\frac{\deg w_{i0}-\deg w_{p0}}{p-i}\right)
\end{equation*}
we conclude that there exists a constant $C_0>0$ such that
\begin{equation}\label{eq:est_2_wnk}
 |W(n,k)|\geq C_0 n^{\deg w_{p0}} k^p\quad\text{for every}\quad k\geq n^{\gamma}\quad\text{and every}\quad n\in\N.
\end{equation}
\smallskip\par
\emph{Step III. Finding a formal solution of (\ref{eq:n}).}

We assume that $f_n(z)=\sum_{k=0}^{\infty}f_{n,k}z^k$. Let $u_n(z)=\sum_{k=0}^{\infty}u_{n,k}z^k$ be a formal solution of (\ref{eq:n}). To calculate $u_{n,k}$ observe that for $s>0$ we have
\begin{multline*}
 P(n,\partial_z)\left(\sum_{k=0}^{\infty}u_{n,k}z^k\right)=\sum_{k=0}^{\infty}W(n,k)u_{n,k}z^k + \sum_{k=0}^{\infty}\sum_{(i,j)\in\tilde{\Lambda}_s\cap\Lambda'}\sum_{l=0}^{\infty}w_{ij}(n)a_{ijl}k^iu_{n,k}z^{k+j+l}\\
 +\sum_{k=0}^{\infty}\sum_{(i,j)\in\Lambda'\setminus\tilde{\Lambda}_s}\sum_{l=0}^{\infty}w_{ij}(n)a_{ijl}k^i u_{n,k}z^{k+j+l}=\sum_{k=0}^{\infty}f_{n,k}z^k.
\end{multline*}

Hence we get the inductive formula on $u_{n,k}$ with respect to $k\in\N_0$.
\begin{multline}
\label{eq:formal}
 W(n,k)u_{n,k}=-\sum_{(i,j)\in\tilde{\Lambda}_s\cap\Lambda'}\sum_{l=0}^{k-j}w_{ij}(n)a_{ijl}(k-j-l)^i u_{n,k-j-l} \\
 -\sum_{(i,j)\in\Lambda'\setminus\tilde{\Lambda}_s}\sum_{l=0}^{k-j}w_{ij}(n)a_{ijl}(k-j-l)^i u_{n,k-j-l} +f_{n,k},
\end{multline}
where we put $u_{n,\tilde{k}}=0$ if $\tilde{k}<0$.
\medskip\par
\emph{Step IV. Estimation of $u_{n,k}$.}

Let
\begin{equation*}
 \beta:=\max_{(i,j)\in\Lambda'}\frac{\deg w_{ij}-\deg w_{p0}}{j}.
\end{equation*}
Observe that $\beta\geq\alpha$. We will show that there exists $B>0$ such that for every sequence of positive numbers $(\tilde{A}(n))_{n\in\N}$ the following implication holds: if
\begin{equation}\label{eq:f_nk}
 |f_{n,k}|\leq \tilde{A}(n)B^k k!^s n^{\beta k}k^{kp}n^{k\deg w_{p0}}\quad\text{for every}\quad k,n\in\N
\end{equation}
then 
\begin{equation}\label{eq:u_nk}
 |u_{n,k}|\leq \frac{3}{C_0}\tilde{A}(n)B^k k!^s n^{\beta k}k^{kp}n^{k\deg w_{p0}}\quad\text{for every}\quad k,n\in\N.
\end{equation}

We will prove (\ref{eq:u_nk}) inductively with respect to $k$ for any fixed $n\in\N$.

By (\ref{eq:snrc}) we see that $|W(n,k)|\geq C_0$. Hence, by (\ref{eq:formal}) and using the estimations (\ref{eq:est_w_ij}) and (\ref{eq:est_a_ijl}) we get
\begin{multline*}
 |u_{n,k}|\leq \sum_{(i,j)\in\tilde{\Lambda}_s\cap\Lambda'}\sum_{l=0}^{k-j}\frac{C_1}{C_0}n^{\deg w_{ij}} A_1 B_1^l l!^s (k-j-l)^i |u_{n,k-j-l}|\\
 +
 \sum_{(i,j)\in\Lambda'\setminus\tilde{\Lambda}_s}\sum_{l=0}^{k-j}\frac{C_1}{C_0}n^{\deg w_{ij}}A_1B_1^l l!^s (k-j-l)^i |u_{n,k-j-l}| +\frac{|f_{kn}|}{C_0}\\
 \leq \frac{A_1 C_1}{C_0} k^p n^{\deg w_{p0}}\Bigg[\sum_{(i,j)\in\tilde{\Lambda}_s\cap\Lambda'}
 \sum_{l=0}^{k-j} n^{\deg w_{ij} - \deg w_{p0}}\frac{(k-j-l)^i l!^s}{k^p}B_1^l|u_{n,k-j-l}|\\
 +
 \sum_{(i,j)\in\Lambda'\setminus\tilde{\Lambda}_s}\sum_{l=0}^{k-j}n^{\deg w_{ij}-\deg w_{p0}}
 \frac{(k-j-l)^il!^s}{k^p}B_1^l|u_{n,k-j-l}|\Bigg]+\frac{|f_{kn}|}{C_0}
\end{multline*}
It means that by (\ref{eq:alpha}) and by the inductive assumption we get
\begin{multline*}
 |u_{n,k}|\leq \frac{|f_{kn}|}{C_0}+\frac{A_1 C_1}{C_0^2} k^p n^{\deg w_{p0}}3\tilde{A}(n)k!^s\\
 \times
 \Bigg[\sum_{(i,j)\in\tilde{\Lambda}_s\cap\Lambda'}\sum_{l=0}^{k-j}n^{j\alpha}
 \frac{(k-j-l)^i (k-j-l)!^s l!^s}{k^p k!^s}B_1^l B^{k-j-l} n^{(k-j-l)\beta}(k-j-l)^{(k-j-l)p}n^{(k-j-l)\deg w_{p0}}\\
 + \sum_{(i,j)\in\Lambda'\setminus\tilde{\Lambda}_s}\sum_{l=0}^{k-j}n^{j\beta}
 \frac{(k-j-l)^i(k-j-l)!^sl!^s}{k^p k!^s}B_1^lB^{k-j-l}n^{(k-j-l)\beta}(k-j-l)^{(k-j-l)p}
 n^{(k-j-l)\deg w_{p0}}
 \Bigg].
\end{multline*}
Since $\beta\geq\alpha$, the first sum in the square brackets multiplied by
$\frac{A_1 C_1}{C_0^2} k^p n^{\deg w_{p0}}3\tilde{A}(n)k!^s$
can be estimated by
\begin{equation*}
\frac{3}{C_0}\tilde{A}(n) B^k k!^s n^{\beta k} k^{kp} n^{k\deg w_{p0}} \frac{A_1 C_1}{C_0} \sum_{(i,j)\in\tilde{\Lambda}_s\cap\Lambda'} \sum_{l=0}^{k-j}B_1^l B^{-j-l}\frac{(k-j-l)^i (k-j-l)!^s l!^s}{k^p k!^s}.
\end{equation*}
Moreover putting 
\begin{equation}\label{eq:B_tilde}
\tilde{B}:=|\Lambda'|
\end{equation}
and using Lemma \ref{le:1} we get
\begin{multline*}
 \frac{A_1 C_1}{C_0} \sum_{(i,j)\in\tilde{\Lambda}_s\cap\Lambda'} \sum_{l=0}^{k-j} B_1^lB^{-j-l} \frac{(k-j-l)^i (k-j-l)!^s l!^s}{k^p k!^s}\\
 \leq
 \frac{A_1 C_1}{C_0}\tilde{B}\bigg(\sum_{j=1}^{\infty}B^{-j}\bigg)\bigg(\sum_{l=0}^{\infty}(B_1/B)^l\bigg)
 \leq\frac{2A_1 C_1}{C_0} \frac{\tilde{B}}{B-1} \leq \frac{1}{3}
\end{multline*}
for 
\begin{equation}\label{eq:B}
B\geq \max(2B_1, 6A_1 C_1 \tilde{B}/C_0 + 1).
\end{equation}

To estimate the second sum in the square brackets multiplied by
$\frac{A_1 C_1}{C_0^2} k^p n^{\deg w_{p0}}3\tilde{A}(n)k!^s$
we put
\begin{equation*}
 s':=\max_{(i,j)\in\Lambda'\setminus\tilde{\Lambda}_s}\frac{i-p}{j}.
\end{equation*}
Of course $s'<s$. Now, we estimate this sum, similarly to the previous one, by
\begin{equation*}
\frac{3}{C_0}\tilde{A}(n) B^k k!^s n^{\beta k}k^{kp}n^{k\deg w_{p0}}\frac{A_1 C_1}{C_0}
\sum_{(i,j)\in\Lambda'\setminus\tilde{\Lambda}_s}\sum_{l=0}^{k-j}B_1^lB^{-j-l}\frac{(k-j-l)^i (k-j-l)!^s l!^s}{k^p k!^s}.
\end{equation*}
Then by Lemma \ref{le:2} we get
\begin{multline*}
 \frac{A_1 C_1}{C_0}
\sum_{(i,j)\in\Lambda'\setminus\tilde{\Lambda}_s}\sum_{l=0}^{k-j}B_1^lB^{-j-l}\frac{(k-j-l)^i (k-j-l)!^s l!^s}{k^p k!^s}\\
\leq \frac{A_1 C_1}{C_0}k^{s'-s} \bigg(\sum_{j=1}^{\infty}B^{-j}\bigg) \bigg(\sum_{l=0}^{\infty}(B_1/B)^l \bigg)
\leq\frac{A_1 C_1}{C_0}k^{s'-s}\frac{2\tilde{B}}{B-1}\leq \frac{1}{3}k^{s'-s}
\end{multline*}
for $B$ satisfying (\ref{eq:B}).

By the above estimations and by (\ref{eq:f_nk}) we receive (\ref{eq:u_nk}).
\medskip\par
\emph{Step V. Better estimation of $u_{n,k}$ for large $k$.}

We assume that $k\geq n^{\gamma}$. Then we may use the estimation of $|W(n,k)|$ given by (\ref{eq:est_2_wnk}).

To improve the estimation calculated in the previous step, we also assume that $k\geq n^{\tilde{\gamma}}$,
where 
\begin{equation*}
 \tilde{\gamma}:=\max_{(i,j)\in\Lambda'}\frac{j(\beta-\alpha)}{s-s'}.
\end{equation*}
Observe that in this case
\begin{equation}\label{eq:gamma_tilde}
 n^{j(\beta-\alpha)}\leq k^{s-s'}\quad\text{for every}\quad (i,j)\in\Lambda'.
\end{equation}
Let $\bar{\gamma}:=\max(\gamma,\tilde{\gamma})$ and $(\tilde{A}(n))_{n\geq 0}$ be a fixed sequence of positive numbers.
We will show that there exists $B>0$ such that for every sequence of positive numbers $(A(n))_{n\in\N}$ satisfying the inequality
\begin{equation*}
 A(n)\geq \frac{3}{C_0}\tilde{A}(n)n^{(\beta-\alpha)n^{\bar{\gamma}}}n^{\bar{\gamma}pn^{\bar{\gamma}}}n^{n^{\bar{\gamma}}\deg w_{p0}}\quad\text{for every}\quad n\in\N,
\end{equation*}
the following implication holds: if
\begin{equation}\label{eq:f_nk_2}
 |f_{n,k}|\leq \tilde{A}(n)B^k k!^s n^{\alpha k}\quad\text{for every}\quad k,n\in\N
\end{equation}
then 
\begin{equation}\label{eq:u_nk_2}
 |u_{n,k}|\leq A(n)B^k k!^s n^{\alpha k}\quad\text{for every}\quad k,n\in\N.
\end{equation}
To prove the implication, we fix $n\in\N$ and we will prove it by the induction with respect to $k$.

If $k\leq n^{\bar{\gamma}}$ then by (\ref{eq:u_nk}) we get
\begin{multline*}
 |u_{n,k}|\leq \frac{3}{C_0}\tilde{A}(n)B^k k!^s n^{\beta k}k^{kp}n^{k\deg w_{p0}}\\\leq
 \frac{3}{C_0}\tilde{A}(n)n^{(\beta-\alpha)n^{\bar{\gamma}}}n^{\bar{\gamma}pn^{\bar{\gamma}}}
 n^{n^{\bar{\gamma}}\deg w_{p0}}
 B^k k!^s n^{\alpha k}\leq A(n)B^k k!^s n^{\alpha k}.
\end{multline*}

Now, we assume that $k\geq n^{\bar{\gamma}}$. 
Similarly to the estimations from the previous step, but using the inequality (\ref{eq:est_2_wnk}) instead of (\ref{eq:snrc}) and the inductive assumption (\ref{eq:u_nk_2}) instead of (\ref{eq:u_nk}), we get
\begin{multline*}
 |u_{n,k}|\leq \frac{|f_{n,k}|}{C_0}+\frac{A_1 C_1}{C_0}
 A(n)k!^s\Bigg[\sum_{(i,j)\in\tilde{\Lambda}_s\cap\Lambda'} \sum_{l=0}^{k-j}
 n^{j\alpha}
 \frac{(k-j-l)^i (k-j-l)!^s l!^s}{k^p k!^s}B_1^lB^{k-j-l}n^{(k-j-l)\alpha}\\
 + \sum_{(i,j)\in\Lambda'\setminus\tilde{\Lambda}_s}\sum_{l=0}^{k-j}n^{j\beta}\frac{(k-j-l)^i(k-j-l)!^sl!^s}{k^p k!^s}B_1^lB^{k-j-l}n^{(k-j-l)\alpha}
 \Bigg].\\
\end{multline*}
So, by (\ref{eq:gamma_tilde}), Lemmas \ref{le:1} and \ref{le:2}, we conclude that
\begin{multline*}
 |u_{n,k}| \leq \frac{|f_{n,k}|}{C_0}\\
 +\frac{A_1 C_1}{C_0}
 A(n)B^kk!^s n^{\alpha k}\Bigg[\sum_{(i,j)\in\tilde{\Lambda}_s\cap\Lambda'}\sum_{l=0}^{k-j}B^{-j}(B_1/B)^l
 +
 \sum_{(i,j)\in\Lambda'\setminus\tilde{\Lambda}_s}\sum_{l=0}^{k-j} n^{j(\beta-\alpha)} k^{s'-s} B^{-j} (B_1/B)^l
 \Bigg]\\
 \leq \frac{|f_{n,k}|}{C_0}+\frac{A_1 C_1}{C_0}
 A(n)B^kk!^s n^{\alpha k}
 \Bigg[\tilde{B}\bigg(\sum_{j=1}^{\infty}B^{-j}\bigg) \bigg(\sum_{l=0}^{\infty}(B_1/B)^l \bigg) + \tilde{B}\bigg(\sum_{j=1}^{\infty}B^{-j}\bigg) \bigg(\sum_{l=0}^{\infty}(B_1/B)^l \bigg)\Bigg]\\ \leq
 \frac{|f_{n,k}|}{C_0}+\frac{A_1 C_1}{C_0}
 A(n)B^kk!^s n^{\alpha k}\frac{4\tilde{B}}{B-1}\leq
 A(n)B^kk!^sn^{\alpha k}
\end{multline*}
for $\tilde{B}$ given by (\ref{eq:B_tilde}) and for $B$ satisfying (\ref{eq:B}).
\medskip\par
\emph{Step VI. Showing that $\alpha>0$ is the smallest possible exponent for which the implication holds.}

Assume that $\alpha$ defined by (\ref{eq:alpha}) is positive and take $(i^*,j^*)\in\tilde{\Lambda}_s\cap\Lambda'$ such that 
\begin{equation*}
\alpha=\frac{\deg w_{i^* j^*}-\deg w_{p0}}{j^*}.
\end{equation*}
Fix $n\in\N$ sufficiently large that $w_{i^* j^*}(\tilde{n})\neq 0$ for every $\tilde{n}\geq n$. In this case there exists $\tilde{D}_1>0$ such that 
\begin{equation}\label{eq:wij*}
|w_{i^* j^*}(n)|\geq\tilde{D}_1 n^{\deg w_{i^* j^*}}.
\end{equation}
We will show that there exist the inhomogeneity
$f_n(z)=\sum_{k=0}^{\infty}f_{n,k}z^k$ with $f_{n,k}$ satisfying (\ref{eq:f_nk_2}) and positive constants $C=C(n)$ and $D$ independent of $n$, such that
\begin{equation}\label{eq:greater}
 |u_{n,k}|\geq C(n)D^kk!^sn^{\alpha k}\quad\text{for every}\quad k=mj^*,\ m\in\N.
\end{equation}
We construct an inhomogeneity $f_n(z)=\sum_{k=0}^{\infty}f_{n,k}z^k$ for which (\ref{eq:greater}) holds inductively with respect to $k$ in the following way: $f_{n,0}=W(n,0)$ and
\begin{equation*}
 f_{n,k}=\sum_{(i,j)\in\Lambda'}
 \sum_{l=0}^{k-j} w_{ij}(n) a_{ijl} (k-j-l)^i u_{n,k-j-l} - w_{i^*j^*}(n) a_{i^*j^*}(0) (k-j^*)^{i^*} u_{n,k-j^*}\quad\text{for}\quad k\in\N,
\end{equation*}
where $u_{n,\tilde{k}}$ are defined by (\ref{eq:formal}) for $\tilde{k}=0,\dots,k-1$ and $u_{n,\tilde{k}}=0$ for $\tilde{k}<0$.

For such constructed inhomogeneity $f_n(z)=\sum_{k=0}^{\infty}f_{n,k}z^k$, by (\ref{eq:formal}) the coefficients $u_{n,k}$ of the formal solution $u_n(z)=\sum_{k=0}^{\infty}u_{n,k}z^k$ of (\ref{eq:n}) satisfy $u_{n,0}=1$ and
\begin{equation}\label{eq:18}
 W(n,k)u_{n,k}=-w_{i^*j^*}(n)a_{i^*j^*}(0)(k-j^*)^{i^*}u_{n,k-j^*}\quad\text{for}\quad k>0,
\end{equation}
where we put $u_{n,\tilde{k}}=0$ for $\tilde{k}<0$.

Since for $k\geq n^{\gamma}$ the leading term of $W(n,k)$ is given by $n^{\deg w_{p0}}k^p$, there exists $\tilde{D}_2>0$ such that 
\begin{equation}\label{eq:wnk_inv}
|W(n,k)|\leq\tilde{D}_2 k^p n^{\deg w_{p0}}\quad\text{for}\quad k\geq n^{\gamma}.
\end{equation}

By (\ref{eq:18}), (\ref{eq:wij*}), (\ref{eq:wnk_inv}) and using Lemma \ref{le:3} we estimate for $k\geq n^{\gamma}$
\begin{multline}
 \frac{|u_{n,k}|}{k!^s}\geq \frac{\tilde{D}_1 |a_{i^*j^*}(0)|}{\tilde{D}_2}\frac{(k-j^*)^{i^*}(k-j^*)!^s}{k^p k!^s}n^{\deg w_{i^*j^*}-\deg w_{p0}}\frac{|u_{n,k-j^*}|}{(k-j^*)!^s}\\\geq\frac{\tilde{D}_1 |a_{i^*j^*}(0)|}{\tilde{D}_2 2^{i^*}}n^{\alpha j^*}\frac{|u_{n,k-j^*}|}{(k-j^*)!^s}.
\end{multline}
Hence there exists $C(n)$ such that for $D=(\frac{\tilde{D}_1 |a_{i^*j^*}(0)|}{\tilde{D}_2 2^{i^*}})^{1/j^*}$ the inequality (\ref{eq:greater}) holds.

From this we conclude that $\alpha>0$ is the smallest possible exponent for which the implication in Theorem \ref{th:1} holds.
Indeed, if the implication would hold for some $\tilde{\alpha}\in(0,\alpha)$ then by the construction of the above inhomogeneity $f_n(z)=\sum_{k=0}^{\infty}f_{n,k}z^k$, this function would satisfy (\ref{eq:f_nk_est}) with $\alpha$ replaced by $\tilde{\alpha}$. Therefore also (\ref{eq:u_nk_est}) should hold with $\alpha$ replaced by $\tilde{\alpha}<\alpha$. But this contradicts (\ref{eq:greater}),
which completes the proof.
\end{proof}

\section{Integro-differential operators in two variables}\label{section:4}

In this section we apply Theorem \ref{th:1} to extend the results
of \cite{lamisu5}, where we study the conditions under which 
the integro-differential operators of the form $P(\partial_t,\partial_z)\partial_t^{-m}$ are automorphisms on the spaces $\C[[t,z]]$ and $\C[[z]]_s[[t]]$. We extend these automorphisms onto the space $(G_s(r(n))_{n\in\N_0}[[t]]$ for some sequence of radii $(r(n))_{n\in\N_0}$ polynomially decreasing to zero as $n$ tends to infinity.

To this end, first we recall the results given in \cite[Section 4]{lamisu5} with a slightly changed notation. We assume that $\Delta$ is a finite subset of indices in
$\N_0\times\N_0$. We consider a partial differential operator defined in \cite[Section 4]{lamisu5} as
\begin{equation}
\label{eq:P}
P(\partial_t,\partial_z)=\sum_{(q,r)\in\Delta}a_{qr}(t,z)\partial_t^q\partial_z^r \quad \textrm{with} \quad a_{qr}(t,z)=\sum_{n=\ord_t(a_{qr})}^\infty a_{q,r,n}(z)t^n\in G_s(r)[[t]].
\end{equation}

We also assume that
\begin{equation*}
m:=\max_{(q,r)\in\Delta}(q-\ord_t(a_{qr}))
\end{equation*}
is greater or equal to zero.

The principal part of the operator $P(\partial_t,\partial_z)$ with respect to $\partial_t$ is given by
\begin{equation*}
P_m(\partial_t,\partial_z):=\sum_{(q,r)\in\Delta_m}a_{q,r,q-m}(z)t^{q-m}\partial_t^q\partial_z^r,
\end{equation*}
where $\Delta_m:=\{(q,r)\in\Delta\colon q-\ord_t(a_{qr})=m\}$.
Then
\begin{equation*}
P_m(\partial_t,\partial_z)\partial_t^{-m}\left(\sum_{n=0}^\infty u_n(z)t^n\right)
=\sum_{n=0}^\infty \tilde{P}_m(n,\partial_z)u_n(z)t^n,
\end{equation*}
where 
\begin{equation}\label{eq:pm}
\tilde{P}_m(n,\partial_z)=\sum_{(q,r)\in\Delta_m}\tilde{a}_{q,r,q-m}(z)z^{\alpha_{q,r}}n(n-1)\cdots(n-(q-m)+1)\partial_z^r
\end{equation}
with $a_{q,r,q-m}(z)=z^{\alpha_{q,r}}\tilde{a}_{q,r,q-m}(z)$ and $\tilde{a}_{q,r,q-m}(0)\neq 0$. Observe that here $\alpha_{q,r}=\ord_z(a_{q,r,q-m}(z))$.

To construct the principal part of the above operator $\tilde{P}_m(n,\partial_z)$ with respect to $\partial_z$, we put
\begin{equation}\label{eq:l}
l:=\min_{(q,r)\in\Delta_m}(\alpha_{q,r}-r)\quad\textrm{and}\quad
\Delta_{m,\tilde{l}}:=\{(q,r)\in\Delta_m\colon \alpha_{q,r}-r=\tilde{l}\}\quad \text{for}\quad \tilde{l}\in\N_0,\quad \tilde{l}\geq l.
\end{equation}
It means that $l$ is the lower ordinate of the Newton polygon 
$N(\tilde{P}_m(n,\partial_z))$ for every $n\in\N_0$.

Then, for given $n\in\N_0$ the principal part of $\tilde{P}_m(n,\partial_z)$ is given by
\begin{equation*}
 \tilde{P}_{m,-l}(n,\partial_z)=\sum_{(q,r)\in\Delta_{m,l}}\tilde{a}_{q,r,q-m}(z)z^{r+l}n(n-1)\cdots(n-(q-m)+1)\partial_z^r,
\end{equation*}
under condition that $\tilde{P}_{m,-l}(n,\partial_z)\neq 0$.

By plugging the power series $u_n(z)=\sum_{k=0}^\infty u_{n,k}z^k$ in $\C[[z]]$ into the previous operator, we get
\begin{equation*}
 \tilde{P}_{m,-l}(n,\partial_z)u_n(z)=\sum_{k=0}^\infty\sum_{(q,r)\in\Delta_{m,l}}\tilde{a}_{q,r,q-m}(z)n(n-1)\cdots(n-(q-m)+1)k(k-1)\cdots(k-r+1)u_{n,k}z^{k+l}.
\end{equation*}

Therefore
\begin{equation*}
 \tilde{P}_{m,-l}(n,\partial_z)u_n(z)=\sum_{k=0}^\infty W_{m,l}(n,k,z)u_{n,k}z^{k+l},
\end{equation*}
where
\begin{equation*}
 W_{m,l}(n,k,z)=\sum_{(q,r)\in\Delta_{m,l}}\tilde{a}_{q,r,q-m}(z)n(n-1)\cdots(n-(q-m)+1)k(k-1)\cdots(k-r+1).
\end{equation*}

Let us recall the results of \cite{lamisu5}.

\begin{theo}[{\cite[Theorem 3]{lamisu5}}]
Any two of the following three conditions entail a third one:
  \begin{enumerate}
 \item[(i)] The operator $P(\partial_t,\partial_z)\partial_t^{-m}$ is a linear automorphism on $\C[[t,z]]$.
  \item[(ii)] The lower ordinate of the Newton polygon $N(\tilde{P}_m(n,\partial_z))$ is equal to $l$ for every $n\in\N_0$, where $l$ is defined by (\ref{eq:l}) and is equal to zero.
  \item[(iii)] The non-resonance condition $W_{m,l}(n,k,0)\ne 0$ holds for every $n,k\in\N_0$.
 \end{enumerate}
\end{theo}

\begin{theo}[{\cite[Theorem 4]{lamisu5}}]\label{th:3}
 Let $s\geq 0$. Assume that the following conditions hold:
 \begin{enumerate}
 \item[(a)] The lower ordinate $l$ of the Newton polygon $N(\tilde{P}_m(n,\partial_z))$ is equal to zero for every $n\in\N_0$.
  \item[(b)] The first positive slope of the Newton polygon $N(\tilde{P}_m(n,\partial_z))$ is greater or equal to $1/s$ for every $n\in\N_0$.
  \item[(c)] The non-resonance condition $W_{m,l}(n,k,0)\ne 0$ holds for every $n,k\in\N_0$.
 \end{enumerate}
 Then the operator $P(\partial_t,\partial_z)\partial_t^{-m}$ is a linear automorphism on $\C[[t,z]]$, which extends to an automorphism on $\C[[z]]_s[[t]]$.
 
 Conversely, if the operator $P(\partial_t,\partial_z)\partial_t^{-m}$ is a linear automorphism on $\C[[t,z]]$, which extends to an automorphism on $\C[[z]]_s[[t]]$, and one of the conditions (a) or (c) is true, then all conditions (a), (b) and (c) are satisfied. 
\end{theo}

Now we are ready to apply Theorem \ref{th:1} to the operator $P(\partial_t,\partial_z)$ defined by (\ref{eq:P}).

Observe that instead of the operator $P(n,\partial_z)$ given in Theorem \ref{th:1} by (\ref{eq:P_n}), we have the operator $\tilde{P}_m(n,\partial_z)$ defined by (\ref{eq:pm}). 

It means that in our case
\begin{equation*}
 p=\max \{r\colon(j,r)\in\Delta_{m,0}\},
\end{equation*}
\begin{equation}\label{eq:deg_wij}
\deg w_{ij}=\max_{(q,i)\in\Delta_{m,j}}(q-m)\quad\text{for}\quad\Delta_{m,j}\neq\emptyset\quad\text{and}\quad \deg w_{ij}=0\quad\text{for}\quad\Delta_{m,j}=\emptyset,
\end{equation}
and
\begin{equation}\label{eq:al}
 \alpha=\max\bigg(0, \max_{j\geq 1,\ js\in\N_0}\frac{\deg w_{p+js,j}-\deg w_{p0}}{j}\bigg)
\end{equation}

We have
\begin{theo}\label{th:4}
 Let $\alpha$ be given by (\ref{eq:al}) and $s\geq 0$. Assume that the following conditions hold:
 \begin{enumerate}
 \item[(a)] The lower ordinate $l$ of the Newton polygon $N(\tilde{P}_m(n,\partial_z))$ is equal to zero for every $n\in\N_0$.
  \item[(b)] The first positive slope of the Newton polygon $N(\tilde{P}_m(n,\partial_z))$ is greater than or equal to $1/s$ for every $n\in\N_0$.
  \item[(c)] There exists $C_0>0$ such that $|W_{m,l}(n,k,0)|\geq C_0$ for every $n,k\in\N_0$ (the strong non-resonance condition).
 \end{enumerate}
 Then the operator $P(\partial_t,\partial_z)\partial_t^{-m}$ is a linear automorphism on $\C[[t,z]]$, which extends to a linear automorphism on $(G_s(a(n+1)^{-\alpha}))_{n\in\N_0}[[t]]$ for some $a>0$.
 Moreover, if $\alpha>0$ then the operator $P(\partial_t,\partial_z)\partial_t^{-m}$ does not extend to a linear automorphism on $(G_s(a(n+1)^{-\tilde{\alpha}}))_{n\in\N_0}[[t]]$ for any $\tilde{\alpha}<\alpha$ and any $a>0$. 
 \end{theo}
\begin{proof}
If all conditions (a), (b), (c) are satisfied then by Theorem \ref{th:3} the operator $P(\partial_t,\partial_z)\partial_t^{-m}$ is a linear automorphism on $\C[[t,z]]$, which extends to an automorphism on $\C[[z]]_s[[t]]$. It means that for every $g=\sum_{n=0}^\infty g_n(z)t^n$ in $\C[[z]]_s[[t]]$  there exists exactly one $u=\sum_{n=0}^\infty u_n(z) t^n$ in $\C[[z]]_s[[t]]$ satisfying $P(\partial_t,\partial_z)\partial_t^{-m}u=g$. By the proof of \cite[Theorem 3]{lamisu5} we conclude that (see \cite[formula (23)]{lamisu5})
\begin{equation}
 \label{eq:u_n}
  \left\{
   \begin{array}{ll}
    \tilde{P}_m(n,\partial_z)u_n(z)=g_n(z)& \textrm{for}\ n=0\\
    \tilde{P}_m(n,\partial_z)u_n(z)=g_n(z)-\sum_{j=1}^{n}\tilde{Q}(n,j,\partial_z)u_{n-j}(z)& \textrm{for}\ n\geq 1.
   \end{array}
  \right.
 \end{equation}
 
 We apply Theorem \ref{th:1} to (\ref{eq:u_n}) with $f_0(z)=g_0(z)$ and 
 \begin{equation}\label{eq:f_n}
 f_n(z)=g_n(z)-\sum_{j=1}^{n}\tilde{Q}(n,j,\partial_z)u_{n-j}(z)\quad \text{for}\quad n\geq 1.
 \end{equation}
 
 We will show inductively that there exists $a>0$ such that if $g\in (G_s(a(n+1)^{-\alpha}))_{n\in\N_0}[[t]]$ then also  $u_n\in G_s(a(n+1)^{-\alpha})$ for every $n\in\N_0$. To this end we put $a=B^{-1}$, where $B$ is a constant depended on the sequence of operators $(\tilde{P}_m(n,\partial_z))_{n\in\N_0}$ from Theorem \ref{th:1}. 
 
  We see that $g_0\in G_s(a)$. Next, applying Theorem \ref{th:1} to the first equation in (\ref{eq:u_n}) with $n=0$ replaced by $n=1$, we conclude that also $u_0\in G_s(a)$.
 
 In the inductive step we assume that $n\geq 1$ is fixed and $u_{n-j}\in G_s(a(n+1)^{-\alpha})$ for $j=1,\dots,n$. Then $\tilde{Q}(n,j,\partial_z)u_{n-j}\in G_s(a(n+1)^{-\alpha})$ for $j=1,\dots,n$. Since also $g_n\in G_s(a(n+1)^{-\alpha})$, by (\ref{eq:f_n}) we see that $f_n\in G_s(a(n+1)^{-\alpha})$. Now, applying Theorem \ref{th:1} to the second equation in (\ref{eq:u_n}) with $n$ replaced by $n+1$, we obtain that also $u_n\in G_s(a(n+1)^{-\alpha})$.
 
 By the induction we conclude that $u_n\in G_s(a(n+1)^{-\alpha})$ for every $n\in\N_0$.  It means that $u(t,z)=\sum_{n=0}^{\infty}u_n(z)t^n\in (G_s(a(n+1)^{-\alpha}))_{n\in\N_0}[[t]]$.
 
 Therefore in this case the operator $P(\partial_t,\partial_z)\partial_t^{-m}$ extends also to a linear automorphism on $(G_s(a(n+1)^{-\alpha}))_{n\in\N}[[t]]$.
 
 To prove the second part of the theorem, suppose, contrary to our claim, that there exists $\tilde{\alpha}\in (0,\alpha)$ such that the operator $P(\partial_t,\partial_z)\partial_t^{-m}$ extends to a linear automorphism on $(G_s(a(n+1)^{-\tilde{\alpha}}))_{n\in\N_0}[[t]]$ for some $a>0$. Then, after repeating the construction of inhomogeneity $f_n(z)$ for $n\in\N_0$ from the proof of Theorem \ref{th:1} (see Step VI), we get the construction of $g=\sum_{n=0}^{\infty}g_n(z)t^n$, which belongs to $(G_s(a(n+1)^{-\tilde{\alpha}}))_{n\in\N_0}[[t]]$, but, as in the proof of Theorem \ref{th:1}, we conclude that then $u=\sum_{n=0}^{\infty}u_n(z)t^n$, satisfying the equation $P(\partial_t,\partial_z)\partial_t^{-m}u=g$, does not belong to $(G_s(a(n+1)^{-\tilde{\alpha}}))_{n\in\N_0}[[t]]$.
Hence, the operator $P(\partial_t,\partial_z)\partial_t^{-m}$ does not extend to a linear automorphism on $(G_s(a(n+1)^{-\tilde{\alpha}}))_{n\in\N_0}[[t]]$ for any $a>0$ and any $\tilde{\alpha}<\alpha$.
\end{proof}

\begin{corol}\label{co:1}
  If the conditions (a)--(c) of Theorem \ref{th:4} hold then the operator $P(\partial_t,\partial_z)\partial_t^{-m}$ extends to a linear automorphism on $G_s[[t]]$
 if and only if $\alpha=0$ or, equivalently if and only if 
\begin{equation*} 
\deg w_{p+js,j}\leq \deg w_{p0}\quad\text{for every}\quad j\in\N_0\quad\text{and}\quad js\in\N_0.
\end{equation*}
\end{corol}
\begin{proof}
It is sufficient to show the first equivalence. The remaining equivalence is obvious.

($\Longleftarrow$) If $\alpha=0$ then by Theorem \ref{th:4} the operator $P(\partial_t,\partial_z)\partial_t^{-m}$ extends to a linear automorphism on  $(G_s(a(n+1)^{0}))_{n\in\N_0}[[t]]=G_s(a)[[t]]$ for some $a>0$, so also on $G_s[[t]]$.

($\Longrightarrow$) If $\alpha>0$ then by Theorem \ref{th:4} with $\tilde{\alpha}=0$ we conclude that the operator $P(\partial_t,\partial_z)\partial_t^{-m}$ does not extend to an automorphism on $(G_s(a(n+1)^{0}))_{n\in\N_0}[[t]]=G_s(a)[[t]]$ for any $a>0$, and therefore the operator $P(\partial_t,\partial_z)\partial_t^{-m}$ does not extend to an automorphism on $G_s[[t]]$.
\end{proof}

To illustrate the above results we give a few examples.
\begin{example}
Let $m\ge0$. We consider the differential operator
$$P(\partial_t,\partial_z)=p_0(t,z)\partial_t^m+p_1(t,z)\partial_t^m\partial_z+p_2(t,z)t\partial_t^{m+1}\partial_z,$$
with 
$$p_0(t,z)=a+zp_{00}(z)+\sum_{n\ge1}p_{0n}(z)t^n\in\C[[t,z]],$$
for some $a>0$, and with $p_{0n}\in\C[[z]]$ for all $n\ge0$,
$$p_1(t,z)=(b+p_{10}(z))z^2+\sum_{n\ge1}p_{1n}(z)t^n\in\C[[t,z]],$$
for some $b>0$, and with $p_{1n}\in\C[[z]]$ for all $n\ge0$, and
$$p_2(t,z)=(c+p_{20}(z))z^h+\sum_{n\ge1}p_{2n}(z)t^n\in\C[[t,z]],$$
with $c>0$ and $p_{2n}\in\C[[z]]$ for all $n\ge0$, and with $h\in\N$, $h>2$. Observe that
$$\tilde{P}_{m}(n,\partial_z)=a+zp_{00}(z)+(b+p_{10}(z))z^2\partial_z+(c+p_{20}(z))z^hn\partial_z,$$
with $N(\tilde{P}_{m}(n,\partial_z))$ being a Newton polygon of one positive slope, and vertices at the points $(0,0)$, $(1,1)$ and $(1,h-1)$, i.e., the first of its positive slopes is 1.
We observe that $W_{m,l=0}(n,k,0)=a$ and therefore the strong non-resonance condition holds. Therefore, the operator $P(\partial_t,\partial_z)\partial_t^{-m}$ is a linear automorphism on $\C[[t,z]]$. 

In our case we have $s=1$, $p=0$, $\deg w_{00}=0$, $\deg w_{11}=0$, $\deg w_{1,h-1}=1$.
Since $\deg w_{11}\leq \deg w_{00}$, we get $\alpha=0$, so by Corollary \ref{co:1} we conclude that the operator $P(\partial_t,\partial_z)\partial_t^{-m}$ extends to an automorphism on $G_s[[t]]$.
\end{example}

\begin{example}\label{ex2}
	Let us now consider an operator of the form
	\begin{equation}
		P(\partial_t,\partial_z)=(\partial_t t)(\partial_z z)-(\partial_t t)^2 z(\partial_z z+1).\label{ex2:operator}
	\end{equation}
	First, notice that $m=0$ and operator (\ref{ex2:operator}) is its own principal part. Moreover, since
	$$
	P(\partial_t,\partial_z)\left(\sum_{n=0}^\infty u_n(z)t^n\right)=\sum_{n=0}^\infty\left[(n+1)(\partial_z z)-(n+1)^2 z(\partial_z z+1)\right]u_n(z)t^n,
	$$ 
	we have  $\tilde{P}_0(n,\partial_z)=(n+1)(\partial_z z)-(n+1)^2 z(\partial_z z+1)$ for every $n\in\N_0$.
	
	The Newton polygon $N(\tilde{P}_0(n,\partial_z))$ for this operator has two vertices, namely $(1,0)$ and $(1,1)$ (see Figure \ref{figure:ex2}), from which we conclude that conditions (a) and (b) from Theorem \ref{th:4} are satisfied with $s=0$.
	
		\begin{figure}[ht]
		\begin{center}
			\includegraphics[scale=0.23]{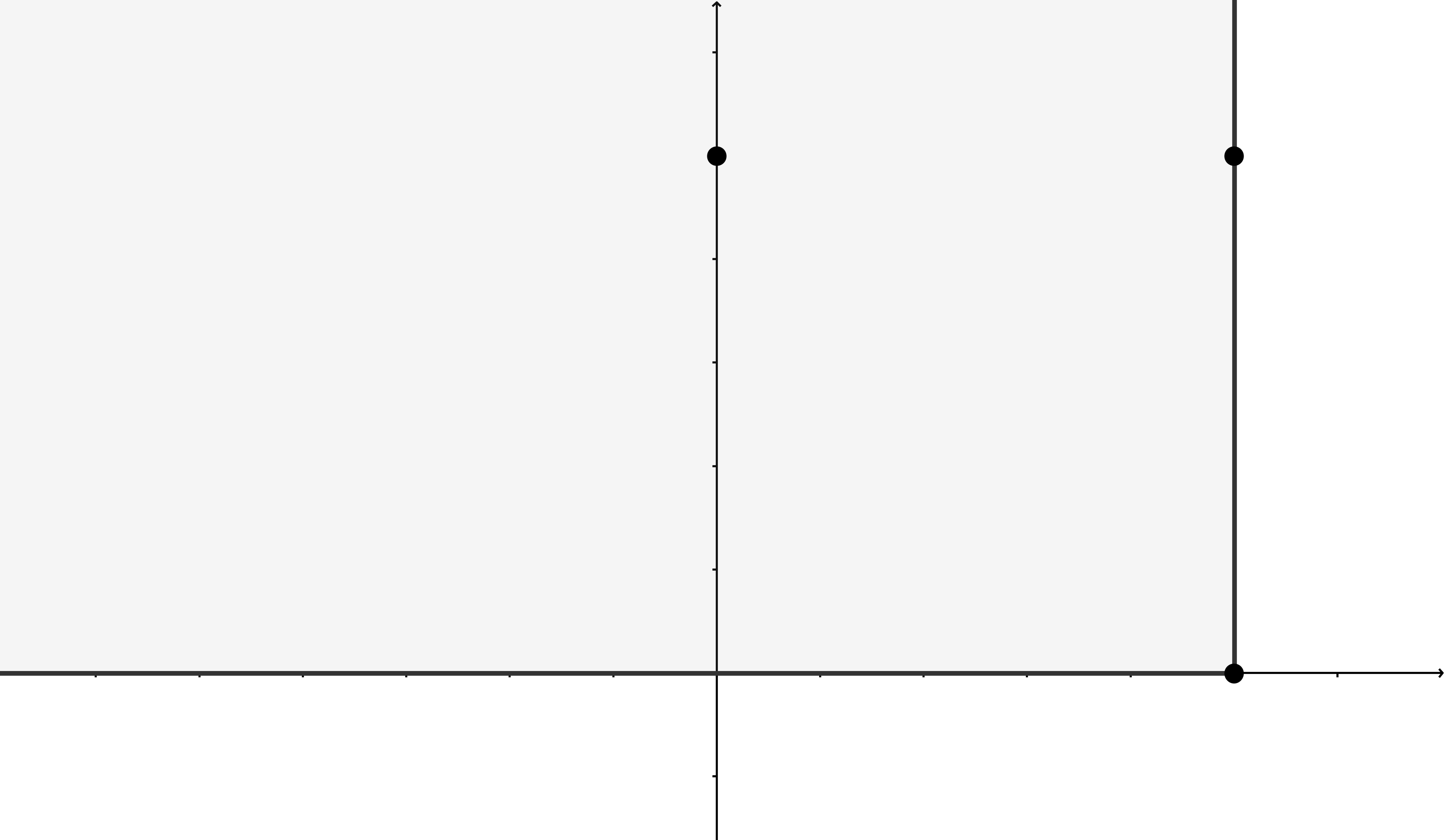}
			\caption{The Newton polygon for $\tilde{P}_0(n,\partial_z)$.\label{figure:ex2}}
		\end{center}
	\end{figure}

 	Moreover, the principal part of $\tilde{P}_0(n,\partial_z)$ is equivalent to
	$$
	\tilde{P}_{0,0}(n,\partial_z)=(n+1)z\partial_z+n+1.
	$$
	Hence, for $u_n(z)=\sum_{k=0}^\infty u_{n,k}z^k$ we receive
	$$
	\tilde{P}_{0,0}(n,\partial_z)(u_n(z))=\sum_{k=0}^\infty (n+1)(k+1)u_{n,k}z^k
	$$
	with $W_{0,0}(n,k,z)=W_{0,0}(n,k,0)=(n+1)(k+1)\ge 1$ for every $n,k\in\N_0$. Hence, condition (c) from Theorem~\ref{th:4} stands. 
	
	 Operator (\ref{ex2:operator}) satisfies the conditions of Theorem \ref{th:4} and, with $p=1$, $w_{11}(n)=-(n+1)^2$ and $w_{10}(n)=n+1$, we receive $\alpha=1$ and can conclude (\ref{ex2:operator}) extends to a linear automorphism on $(G_0(a(n+1)^{-1}))_{n\in\N_0}[[t]]$, but does not extend to a linear automorphism on $(G_0(a(n+1)^{-\tilde{\alpha}}))_{n\in\N_0}[[t]]$ for any $\tilde{\alpha}<1$.
	 \par
	 We will also show this last conclusion directly.
	 
	 This is indeed the case. To show this fact, let us consider an equation of the form
	\begin{equation*}
		P(\partial_t,\partial_z)u(t,z)=f(t,z)
	\end{equation*}
	with $f(t,z)$ given as a formal power series $f(t,z)=\sum_{n=0}^\infty f_n(z)t^n=\sum_{n=0}^\infty\sum_{k=0}^\infty f_{n,k}t^nz^k$. Then the coefficient of its formal power series solution $u(t,z)=\sum_{n=0}^\infty u_n(z)t^n=\sum_{n=0}^\infty\sum_{k=0}^\infty u_{n,k}t^nz^k$ are given by the recurrence formula:
	$$
	(n+1)(\partial_z z)u_n(z)-(n+1)^2 z(\partial_z z+1)u_n(z)=f_n(z)\ \textrm{for }n\in\N_0,
	$$ 
	which for any fixed $n\in\N_0$ gives us
	$$
	(n+1)u_{n,0}=f_{n,0}
	$$
	and
	$$
	(n+1)(k+1)u_{n,k}-(n+1)^2u_{n,k-1}=f_{n,k}\textrm{ for }k\ge 1.
	$$

	Let us now put $f_{n,0}=n+1$ and $f_{n,k}=0$ for $k\ge 1$. Then we receive
	$$
	\begin{cases}
		u_{n,0}=1\\
		u_{n,k}=(n+1)u_{n,k-1}\textrm{ for }k\ge 1
	\end{cases}
	$$
	for every $n\in\N$. Hence, 
	\begin{equation*}
	u_n(z)=\sum_{k=0}^\infty(n+1)^k z^k=\frac{1}{1-(n+1)z}\in\Oo(D_{1/(n+1)})\quad \text{for every}\quad n\in\N.
	\end{equation*}
	We see that $u_n(z)\in\Oo(D_{1/(n+1)})$ and $u_n(z)$ does not extend analytically to any disc $D_{r(n)}$, where $r(n)>1/(n+1)$.
    \par
    It means that 
    \begin{equation}\label{eq:u_example}
    u(t,z)=\sum_{n=0}^{\infty}u_n(z)t^n=\sum_{n=0}^{\infty}\frac{t^n}{1-(n+1)z}\in(G_0(a(n+1)^{-1}))_{n\in\N_0}[[t]],
    \end{equation}
    but
    we cannot replace the exponent $-1$ by a greater one.
    \par
    We have shown that for the inhomogeneity $f(t,z)=\sum_{n=0}^{\infty}(n+1)t^n$, which belongs to the space $(G_0(a(n+1)^{-\tilde{\alpha}}))_{n\in\N_0}[[t]]$ for any $\tilde{\alpha}\in\R$ and any $a>0$,
    a unique formal power series solution $u(t,z)$ of the equation $P(\partial_t,\partial_z)u=f$ does not belong to the space $(G_0(a(n+1)^{-\tilde{\alpha}}))_{n\in\N_0}[[t]]$ for any $\tilde{\alpha}<1$ and any $a>0$.
    Therefore, the operator
    $$P(\partial_t,\partial_z)\colon (G_0(r(n+1)^{-\tilde{\alpha}}))_{n\in\N_0}[[t]] \longrightarrow (G_0(a(n+1)^{-\tilde{\alpha}}))_{n\in\N_0}[[t]]$$
    is not surjective for any $\tilde{\alpha}<1$ and any $a>0$, since 
    $$f(t,z)\not\in \text{Image}\Big(P(\partial_t,\partial_z), (G_0(a(n+1)^{-\tilde{\alpha}}))_{n\in\N_0}[[t]]\Big).$$
\end{example}

\begin{example}\label{ex3}
	Let us now consider an operator
	\begin{equation}
	P(\partial_t,\partial_z)=(\partial_t t)(\partial_z z)-(\partial_t t)^\mu z^\nu (\partial_z z+\nu)\label{ex3:operator}
	\end{equation}
	where $\mu,\nu\in\N$ and $\mu\ge 2$.
	
	This case for $\mu=2$ and $\nu=1$ becomes the one presented in Example \ref{ex2} for the operator~ (\ref{ex2:operator}). Moreover, we have
	$$
	P(\partial_t,\partial_z)\left(\sum_{n=0}^\infty u_n(z) t^n\right)=\sum_{n=0}^\infty\left[(n+1)(\partial_z z)-(n+1)^\mu z^\nu(\partial_z z+\nu)\right]u_n(z)t^n.
	$$
	Hence, $\tilde{P}_0(n,\partial_z)=(n+1)(\partial_z z)-(n+1)^\mu z^\nu(\partial_z z+\nu)$ for every $n\in\N_0$, with its Newton polygon determined by points $(1,0)$ and $(1,\nu)$ corresponding to terms $\partial_z z$ and $z^\nu(\partial_z z)$ respectively. The shape of the Newton polygon is then similar to the one from Example \ref{ex2} with its two vertices being $(1,0)$ and $(1,\nu)$, and we again have $s=0$. It is also easy to notice that $W_{0,0}(n,k,0)=(n+1)(k+1)\ge 1$ for every $n,k\in\N_0$.

	With $p=1$, $w_{10}(n)=n+1$ and $w_{1\nu}=-(n+1)^\mu$ we receive $\alpha=\frac{\mu-1}{\nu}$ and from Theorem~\ref{th:4} it follows that the operator~(\ref{ex3:operator}) can be extended to a linear automorphism on $(G_0(a(n+1)^{-(\mu-1)/\nu}))_{n\in\N_0}[[t]]$ for a certain constant $a>0$, but this operator does not extend to a linear automorphism on $(G_0(a(n+1)^{-\tilde{\alpha}}))_{n\in\N}[[t]]$ for any $\tilde{\alpha}<(\mu-1)/\nu$ and any $a>0$. 
	
	We will show this last conclusion directly.
	Let $u(t,z)=\sum_{n=0}^{\infty}u_n(z)t^n=\sum_{n=0}^\infty\sum_{k=0}^\infty u_{n,k}t^n z^k$ be a formal power series solution of the equation:
	\begin{equation*}
		P(\partial_t,\partial_z)u=f(t,z),
	\end{equation*}
	with the inhomogeneity $f(t,z)$ also given by a formal power series $f(t,z)=\sum_{n=0}^{\infty}f_n(z)=\sum_{n=0}^\infty\sum_{k=0}^\infty f_{n,k}t^n z^k\in\Oo(D_1)[[t]]$. Then, for any fixed $n$, the coefficients $u_n(z)$ are given by a recursive formula:
	$$
	(n+1)(\partial_z z)u_n(z)-(n+1)^\mu z^\nu(\partial_z z+\nu)u_n(z)=f_n(z).
	$$
	From this, for any fixed $k\ge \nu$ we receive:
	$$
	u_{n,k}=\frac{f_{n,k}}{(n+1)(k+1)}+(n+1)^{\mu-1}u_{n,k-\nu}.
	$$
	Moreover, we receive $u_{n,k}=\frac{f_{n,k}}{(n+1)(k+1)}$ for $k=1,2,\ldots,\nu-1$. For the sake of simplicity, let us fix $f_{n,0}=n+1$ and $f_{n,k}=0$ for $k\ge 1$. Then $u_{n,0}=1$ and 
	$$
	u_{n,\nu k}=(n+1)^{(\mu-1)k}\ \textrm{ for }k\ge 1.
	$$
	Hence,
	$$
	u_n(z)=\sum_{k=0}^\infty (n+1)^{(\mu-1)k}z^{\nu k},
	$$
	from which it follows for every $n\in\N_0$ that $u_n(z)\not\in\Oo(D_{r(n)})$ with $r(n)=a(n+1)^{-\tilde{a}}$ for any $\tilde{\alpha}<(\mu-1)/\nu$ and any $a>0$. This agrees with the conclusion following from Theorem~\ref{th:4}.
\end{example}
\section{Final remarks}\label{section:final_remarks}
The presented results shed new light on the appearance of formal solutions with coefficients being holomorphic (or, more generally, Gevrey) functions on shrinking discs described by a sequence of positive radii $(r(n))_{n\geq 0}$ decreasing to zero. In case of partial differential equations we get here sequences of radii $r(n)=a(1+n)^{-\alpha}$, $n\in\N_0$, that polynomially decrease to zero, with the exponent $\alpha$, which is given by (\ref{eq:alpha}), being dependent on the relation between the derivatives $(t\partial_t)$ and $(z\partial_z)$ in the principal part of the operator $P(\partial_t,\partial_z)\partial_t^{-m}$ with respect to both variables.

In the case of $q$-difference equations for $q>1$ (see \cite{cala,tahara3,lama4}) similar phenomena have completely different character and are caused by $q$-difference operator which decreases $q$ times a radius of disc on which a given function is holomorphic.

In the future we are planning to describe such phenomena for moment partial differential operators $P(\partial_{m_1,t},\partial_{m_2,z})$ where sequences of radii $(r(n))_{n\geq 0}$ of shrinking discs depend both on the relation
between $(t\partial_{m_1,t})$ and $(z\partial_{m_2,z})$, and on a decreasing of radius of holomorphicity by the moment derivative $\partial_{m_2,z}$, as in the case of $q$-difference operator with respect to $z$.

\section{Auxiliary lemmas}\label{section:lemmas}
This section contains technical lemmas needed in the proof of Theorem~\ref{th:1}.
\begin{lemma}\label{le:1}
	If $s\geq 0$, $i,p\in\N_0$, $j\in\N$ and $js\geq i-p$ then
	\begin{equation}\label{eq:le_1}
		\frac{(k-j-l)^i (k-j-l)!^s l!^s}{k^p k!^s}\leq 1 
	\end{equation}
	for every $k\in\N$ and $l\in\N_0$ satisfying $k-j-l\geq 0$.
\end{lemma}
\begin{proof}
	If $k-j-l=0$ then (\ref{eq:le_1}) is obvious.
	So, we assume that $k-j-l>0$. In this case we have
	\begin{multline*}
		\frac{(k-j-l)^i(k-j-l)!^sl!^s}{k^p k!^s}\leq \frac{(k-j)^i}{(k-j)^p}\frac{(k-j)!^s}{k!^s}\leq (k-j)^{i-p}\bigg(\frac{(k-j)!}{k!}\bigg)^s\\
		\leq
		\bigg(\frac{k!}{(k-j)!}\bigg)^{\frac{i-p}{j}}
		\bigg(\frac{(k-j)!}{k!}\bigg)^{\frac{i-p}{j}}\leq 1.
	\end{multline*}
\end{proof}

\begin{lemma}\label{le:2}
	Assume that $s\geq 0$, $i,p\in\N_0$, $j\in\N$ and 
	there exists $s'<s$ such that  $js'\geq i-p$. Then
	\begin{equation}\label{eq:le_2}
		\frac{(k-j-l)^i(k-j-l)!^sl!^s}{k^p k!^s}\leq k^{s'-s}
	\end{equation}
	for every $k\in\N$ and $l\in\N_0$ satisfying $k-j-l\geq 0$.
\end{lemma}
\begin{proof}
	If $k-j-l=0$ then (\ref{eq:le_2}) is trivial.
	So, we may assume that $k-j-l>0$. Then we get
	\begin{multline*}
		\frac{(k-j-l)^i(k-j-l)!^sl!^s}{k^p k!^s}\leq 
		\bigg(\frac{k!}{(k-j)!}\bigg)^{\frac{i-p}{j}}
		\bigg(\frac{(k-j)!}{k!}\bigg)^{s}\leq
		\bigg(\frac{k!}{(k-j)!}\bigg)^{s'}
		\bigg(\frac{(k-j)!}{k!}\bigg)^{s}\\
		\leq\bigg(\frac{(k-j)!}{k!}\bigg)^{s-s'}\leq \frac{1}{k^{s-s'}}.
	\end{multline*}
\end{proof}

\begin{lemma}\label{le:3}
	If $s\geq 0$, $i,p\in\N_0$, $j\in\N$, and $js= i-p$ then
	\begin{equation}\label{eq:le_3}
		\frac{(k-j)^i(k-j)!^s}{k^p k!^s}\geq \frac{1}{2^i} 
	\end{equation}
	for every $k\in\N$ satisfying $k\geq 2j$.
\end{lemma}
\begin{proof}
	If $k\geq 2j$ then $2(k-j)\geq k$. Hence
	\begin{equation*}
		\frac{(k-j)^i}{k^p}\geq \frac{1}{2^i}k^{i-p}\geq \frac{1}{2^i}\bigg(\frac{k!}{(k-j)!}\bigg)^{\frac{i-p}{j}}.
	\end{equation*}
	It means that
	\begin{equation*}
		\frac{(k-j)^i(k-j)!^s}{k^p k!^s}\geq \frac{1}{2^i} \bigg(\frac{k!}{(k-j)!}\bigg)^{\frac{i-p}{j}}\bigg(\frac{(k-j)!}{k!}\bigg)^{\frac{i-p}{j}}\geq \frac{1}{2^i}.
	\end{equation*}
\end{proof}

\section*{Acknowledgements} The authors wish to express their gratitude to professor Masafumi Yoshino for very helpful comments concerning the non-resonance conditions.

\end{document}